\documentclass[preprint,12pt]{elsarticle}

\usepackage{amssymb}
\usepackage{amsmath}

\usepackage{amssymb}
\usepackage{graphicx}
\usepackage{dcpic,pictex}
\usepackage{tikz}
\usetikzlibrary{graphs,graphs.standard,quotes}
\usetikzlibrary{arrows}

\newcommand{\C}{{\mathbb C}}
\newcommand{\rank}{{\rm rank}}

\newcommand{\conv}{{\rm conv}}

\newtheorem{thm}{Theorem}[section]

\newtheorem{rmk}{Remark}[section]
\newtheorem{corol}{Corollary}[section]
\newtheorem{lem}{Lemma}[section]

\newtheorem{Exa}{Example}[section]

\newproof{pf}{Proof}

\numberwithin{equation}{section}

\numberwithin{figure}{section}

\journal{Journal of Mathematical Analysis and Applications}

\begin{document}

\begin{frontmatter}



\title{On the ellipticity of the higher rank numerical range}



\author[1]{N. Bebiano%
\fnref{fn1}}
\ead{bebiano@mat.uc.pt}
\affiliation[1]{organization={CMUC, Departament of Mathematics},
addressline={University of Coimbra},
postcode={3001-501},
city={Coimbra},
country={Portugal}}

\author[2]{R. Lemos\corref{cor1}%
\fnref{fn2}}
\ead{rute@ua.pt}
\affiliation[2]{organization={Center for Research and Development in Ma\-the\-ma\-tics and Applications (CIDMA), Department of Mathematics},
            addressline={University of Aveiro}, 
postcode={3810-193},
city={Aveiro}, 
            country={Portugal}}

\author[3]{G. Soares%
\fnref{fn3}}
\ead{gsoares@utad.pt}
\affiliation[3]{organization={CMAT-UTAD, Department of  Mathematics},
            addressline={University of
Tr\'as-os-Montes  e Alto Douro}, 
postcode={5000-911},
            city={Vila Real},             
            country={Portugal}}

\cortext[cor1]{Corresponding author}
\fntext[fn1]{The author was partially supported by the Centre for Mathe\-matics of the
University of Coimbra, funded by the Portuguese Government through FCT/MCTES, within the Project UID/00324/2025 (https://doi.org/10.54499/UID/00324/2025).
}
\fntext[fn2]{The author was supported by CIDMA (https://ror.org/05pm2mw36)
under the Portuguese Foundation for Science and Technology 
(FCT, https://ror.org/00snfqn58), Grants UID/04106/2025 (https://doi.org/10.54499/UID/04106/2025) and UID/PRR/04106/2025.}
\fntext[fn3]{The author was partially financed by Portuguese Funds through FCT (Fundação para a Ciência e a Tecnologia), within the Project UID/00013/2025 (https://doi.org/10.54499/UID/00013/2025).}

\begin{abstract}
 The higher rank numerical range is a concept that generalizes the classical numerical range
 and it has application in  quantum error correction. We investigate these sets for $2$-by-$2$ block matrices
 with associated Kippenhahn curves consisting of ellipses (and eventually points).
 As a consequence, elliptical higher rank numerical range results are derived in a unified way, using an approach developed by Spitkovsky \emph{et al}.     	 
\end{abstract}




\end{frontmatter}

\section{Introduction} \label{Sec1}

 Let $M_{m,n}$ be the set of
$m\times n$ complex matrices, with 
$M_{n,n}$ abbreviated to $M_n$.
 
 Let $B({\mathcal H})$ 
 be the algebra of bounded linear operators on a Hilbert space $\mathcal H$, identified with $M_n$ if $\mathcal H$ is $n$-dimensional. 
The classical {\it numerical range}  of $A\in B({\mathcal H})$ is  defined and denoted as
$$
W(A) :=\, \left\{\langle A\,x, x\rangle :\ x\in{\mathcal H}\, , \langle x,x \rangle =1\right\}.
$$
This concept was  introduced in the second decade of the $20$th century   by    Toeplitz \cite{T1918} and Hausdorff   \cite{H1919}, and it has been extensively investigated by pure and applied scientists. 
The numerical range has attractive properties (see e.g.\! \cite{GR}), 
such as  convexity, asserted by the Toeplitz-Hausdorf theorem \cite{H1919,T1918}, and the closure of $W(A)$ containing the spectrum of $A$. 
Some  applications in Physics were considered in \cite{BdP1998, BLdP2004}.

Motivated by problems in quantum error correction, Choi, Kribs and \.{Z}yczkowski \cite{CKZ2006, CKZ2006A} introduced  for $A\in M_n$ and any integer $k=1, \ldots, n$, the  {\it rank-k} or {\it higher rank numerical range} of $A$. For $A\in B({\mathcal H})$ and $k \leq \dim \mathcal H$, this set is defined and denoted as
$$\Lambda_k(A) := \{\lambda\in {\mathbb C}: PAP =\lambda P \mbox{ for some rank-}k \mbox{ orthogonal projection } P\in B({\mathcal H} \}.$$
An equivalent definition for $A \in M_n$ is
\begin{equation}\label{equidef}
\Lambda_k(A) \,=\, \{\lambda\in{\mathbb C:\, X^*}AX =\lambda I_k \mbox{ for some } X\in M_{n,k} \mbox{ such that } X^*X=I_k \}.
\end{equation}

From the previous characterization, it is clear that the rank-$1$ numerical range $\Lambda_1(A)$ coincides with $W(A)$ and that the higher rank numerical ranges of $A\in M_n$ form a decreasing sequence of compact sets: 
$$\Lambda_n(A)\,\subseteq\,\cdots\,\subseteq\,\Lambda_1(A)\,=\,W(A).$$
In particular, $\Lambda_n(A)\neq \emptyset$ if and only if $A=\lambda I_n$, in which case all the sets $\Lambda_k(A)=\{\lambda\}$.
For $k>n/2$, $\Lambda_k(A)$  is either the empty set or a singleton $\{\lambda_0\}$, in the latter case $\lambda_0$ being an eigenvalue of $A$ with geometric multiplicity at least $2k-n$ \cite[Proposition 2.2]{CKZ2006}. 
  Li, Poon and Sze \cite{Li2009} showed that $\Lambda_k(A)$ is non-empty if
 $k <n/3+1$. 

As usual, for $A\in M_n$
we write $A=\Re(A)+i\Im(A)$, where
 $$
  \Re(A):=(A+A^*)/2\qquad \hbox{and}\qquad \Im (A):=(A-A^*)/ 2i
 $$
are the {\it real} and {\it imaginary parts} of $A$, respectively. 
We denote by $\lambda_1(\theta), \ldots, \lambda_n(\theta)$ the eigenvalues, counting the multiplicities, of the Hermitian matrix 
 $$ 
  \Re(e^{-i\theta} A) \,=\, \Re(A)\cos\theta+\Im(A)\sin\theta,\qquad \theta\in[0,2\pi).
 $$

The convexity of higher rank numerical ranges holds for $k=1$ and $k>n/2$, and the convexity for intermediate values of $k$ was proved by Woerdeman in \cite{W2008}, and independently by Li and Sze  in \cite{Li2008}. In this paper, the authors obtained
 the following description of the higher rank numerical range of $A\in M_n$ as the intersection of half-planes:
\begin{equation}\label{i}
\Lambda_k(A) \,= \bigcap_{\theta\in [0, 2\pi)} \left\{z\in \mathbb{C}: \ \Re(e^{-i\theta} z )\leq \lambda_k{(\theta)}\right\},
\end{equation}
where $\lambda_k(\theta)$ 
denotes  the $k$-th largest  eigenvalue, counting the multiplicities, of the Hermitian matrix $\Re(e^{-i\theta}A)$, for each $\theta\in [0,2\pi)$, and $k=1, \dots, n$. 
For $A\in M_n$ normal with eigenvalues $\lambda_1,\ldots,\lambda_n$, then (\ref{i}) yields 
 $$
 \Lambda_k(A) \,= \bigcap_{1\leq j_1 \leq  \cdots\leq j_{n+k-1}\leq n} \conv\,\{\lambda_{j_1},\dots,\lambda_{j_{n-k+1}}\}\,,\
 $$ 
 where $\conv\, S$ denotes the convex hull of the set $S$.
For $A\in M_n$ Hermitian with eigenvalues $\lambda_1 \geq\dots\geq\lambda_n$, then
 $$
  \Lambda_k(A) \,=\,\big[\lambda_{n-k+1}, \lambda_k\big]
 $$
when $\lambda_{n-k+1}<\lambda_k$, $\Lambda_k(A)$ is a singleton if $\lambda_k=\lambda_{n-k+1}$ and an empty set, otherwise.
 
Kippenhahn \cite{K1951} proved that $W(A)$ is the convex hull of a certain algebraic curve $C(A)$,
associated with  the matrix $A$. 
 Not only $W(A)$ but also any $\Lambda_k(A)$ can be described in terms of $C(A)$  (see, e.g.\! \cite{WG2013, BPS2022, JC2026}).
 The boundary lines of the half-planes in the right hand side of the representation \eqref{i}, that is,
 $$\ell_{\theta,k}\,=\,\big\{z\in {\mathbb C}:\ \Re(e^{-i\theta}z)=\lambda_k(\theta)\big\},\qquad \theta \in [0, 2\pi), $$
when taken for all $k=1, \ldots, n,$ form a family,  the envelope of which may determine the boundary generating curve $C(A)$ of the  numerical range of $A$. 

In general, it is difficult to  characterize $\Lambda_k(A)$ for any $k$ and, in particular, to describe  the boundary of $\Lambda_k(A).$  

We  mainly focus on complex matrices of the form
\begin{equation}\label{bloco}
A=\left[\begin{array}{cc}
\alpha I_{r}  & C\\
D & \beta I_{n-r}
\end{array}\right],\qquad 0<r<n.
\end{equation}
Clearly, $\Lambda_k(A)\neq \emptyset$ for $k \leq \max\{r, n-r\}$.
 Our goal is to investigate the higher rank numerical range  of 
matrices of this form \eqref{bloco}.

The paper is organized as follows. Section \ref{Sec2} contains useful preliminary results on block matrices of type \eqref{bloco}. 
In Section \ref{Sec3},   the ellipticity of  higher rank numerical ranges 
of  block  matrices of type \eqref{bloco}, such that  $\Re(e^{-i\theta}A)$ happens to be
 unitarily reducible to the direct sum of $2$-by-$2$ (and eventually $1$-by-$1$) matrices, is investigated. In such cases,
 $C(A)$ splits into at most $n/2$ elliptical components, each solely responsible for the respective higher rank numerical range. 
Using the obtained results, the ellipticity of the higher rank numerical ranges of several classes of matrices is derived in an unified way.
\smallskip

\section{Preliminaries}\label{Sec2}

\smallskip

\noindent 
   For
$A \in M_n$ and $1\leq k\leq n$, the following elementary properties of $\Lambda_k(A)$ 
can easily be checked. 

\smallskip

\begin{enumerate}
	\item[P1.] $\Lambda_k(\alpha A+\beta I_n) =\alpha \Lambda_k(A)+\beta$ for any $\alpha, \beta\in \mathbb{C}$ ({\it translational property}).\smallskip
\item[P2.]   $\Lambda_k(A)$ is {\it unitarily invariant}: 	$\Lambda_k(U^*AU)=\Lambda_k(A)$  for any unitary matrix $U\in M_n$.\smallskip
\item[P3.]   $\Lambda_k(A_1\oplus A_2)\supseteq  {\rm conv} \left( \Lambda_k(A_1)\cup \Lambda_k(A_2)\right)$ for any $A_1, A_2\in M_n$. 
\smallskip
\item[P4.]   $\Lambda_{k_1+k_2}(A_1\oplus A_2) \supseteq \Lambda_{k_1}(A_1)\cap \Lambda_{k_2}(A_2)$ for any $k_1,k_2\in \{1, \ldots, n\}$.
\end{enumerate}

\medskip

We will concentrate on matrices of the form \eqref{bloco}.
By property $P_1,$ when investigating $\Lambda_k(A)$ it suffices to consider matrices with $\alpha+\beta=0$.
Without loss of generality, we can moreover assume $n\leq 2r$, because if $n>2r$ by an adequate unitary similarity of $A$, which by property $P_2$ preserves  $\Lambda_k(A)$, we may switch the  matrix to the following one:
\begin{equation}\label{bloco1}
\left[
\begin{array}{cc}
\beta I_{n-r}  & D\\
C & \alpha I_{r}
\end{array}\right].
\end{equation}
 All the matrices $A\in M_2$ have this form, and the  {\it Elliptical Range Theorem}  (see e.g.\! \cite{Li1996})
 states that for $A\in M_2$ with eigenvalues $\lambda_1,\lambda_2,$
 $W(A)$ is an elliptical disc centered at $\frac{1}{2}\mbox{Tr}\,A,$
 with foci at  the eigenvalues,   
   major and minor  axes of lengths, respectively
$$
\left(\mbox{Tr}(A^*A) -2\Re(\lambda_1\bar\lambda_2)\right)^{\frac{1}{2}}\qquad \hbox{and}\qquad 
\left(\mbox{Tr}(A^*A)-|\lambda_1|^2-|\lambda_2|^2\right)^{\frac{1}{2}}.
$$
In this case, $W(A)$ degenerates into the line segment joining $\lambda_1$ and $\lambda_2$  if and only if $A$ is normal, and 
it  reduces to a singleton if and only if $A$ is a scalar matrix. 
The elliptical shape of the $2$-by-$2$ case persists for the numerical range of certain classes of matrices independently of the matrix size (see e.g.\!\!  \cite{BS2004, CH2012, GS2021}). Other structured matrix classes \cite{BLS2023b,BLS2024, CS2015} have different geometric behaviour when the order increases and several matrices whose numerical range is the convex hull of ellipses were investigated (see e.g.\!\! \cite{BPN2013, BPSK2021, JS2024}).

As already observed, the eigenvalues of the Hermitian matrix 
 $$H_\theta(A):=\,\Re(e^{-i\theta}A), \quad \theta\in [0,2\pi),$$ 
play a central role in our study.
 For matrices of type \eqref{bloco}  they have been  obtained in~\cite{GS2021} in terms of those of the matrix 
\begin{equation}\label{Mtheta}
	M_{C,D}(\theta) :=\,
	C^*C+D\,D^*+\,2\,\Re(e^{-2i\theta}DC),  \qquad  \theta\in[0,2\pi),
\end{equation}
if $n\leq 2r$. If $n> 2r$, switching to the matrix (\ref{bloco1}),  they may be derived from the eigenvalues of  
$M_{D,C}(\theta)$.
We observe that $M_{C,D}(\theta)$ and $M_{D,C}(\theta)$ share the same non-zero eigenvalues.
The proof here presented exhibits the unitary similarity of $H_\theta(B)$, where
\begin{equation}\label{B}
  B\,=\,A -\tfrac{1}{2}(\alpha+\beta)\,I_n,
\end{equation}
to a direct sum of small sized  blocks, which is crucial to derive our results.
For simplicity of notation, we write $M(\theta)=M_{C,D}(\theta)$ if $n\leq 2r$ and $M(\theta)=M_{D,C}(\theta)$ if $n> 2r$,
\begin{equation}\label{w}
 w\,=\,\tfrac{1}{2}(\alpha-\beta) \qquad 
\hbox{and} \qquad w_\theta\, =\,\Re(w\,e^{-i\theta}).
\end{equation}

\medskip

\begin{lem}\label{eigs}  
Let $A$ be of the form \eqref{bloco} 
and $\theta\in [0,2\pi)$. The following holds.

\begin{itemize}

\item[\bf (i)] The eigenvalues of $M(\theta)$ are  
	$\mu_j(\theta)=4\,s_j^2(\theta)$,  $1\leq j\leq p$,  and  $\mu_j(\theta)=0$, $j>p$, with  $s_1(\theta),\dots,s_{p}(\theta)$  the non-zero singular values, counting the multiplicities, of 
	$$N_{\theta}=\tfrac{1}{2}\left(e^{-i\theta}C+e^{i\theta}\,D^*\right).$$ 
\item[\bf (ii)] For $B$ given by {\rm (\ref{B})}, the matrix 
$H_\theta(B)$ is unitarily similar to the direct sum 
 $$ 
 {B}_1\,\oplus\,\cdots\,\oplus\,{B}_p\,\oplus\,w_\theta\,I_{r-p}\,\oplus\,(-w_\theta)\,I_{n-r-p}\,,
 $$ 
 where the non-scalar blocks have the form	
  \begin{eqnarray*}
	B_j\,= \left[\begin{matrix}w_\theta & s_j(\theta)\\
	s_j(\theta) & -w_\theta\\
	\end{matrix} \right], \qquad j=1, \ldots,p.
	\end{eqnarray*}	
\item[\bf (iii)]  If $n< 2r$  ($n>2r$), then the  eigenvalues of $H_\theta(A)$  are $\Re (e^{-i\theta}\alpha)$ 
(resp.\! $\Re (e^{-i\theta}\beta)$)
and
	\begin{equation}\label{lambdaj}
	\tfrac{1}{2}{\Re \big(e^{-i\theta}}(\alpha+\beta)\big) \pm\tfrac{1}{2}\sqrt{4\,w_{\theta}^2 + \mu_j(\theta)} \,,
	\end{equation} 
for $j=1,\ldots, n-r$  (resp.\! $j=1,\ldots,r$).	
If $n=2r$,  the eigenvalues of $H_\theta (A)$   are those 	of the form {\rm (\ref{lambdaj})} for  $j=1,\ldots,n/2$. 
	\end{itemize}
\end{lem}

\begin{pf} {\bf (i)} If $n\leq 2r$, the result is an obvious consequence of the fact $M(\theta)=4\,N^*_\theta N_{\theta}$. Otherwise, interchanging $C$ and $D$, the conclusion follows.
\smallskip

\noindent {\bf (ii)} 
	 Letting $A=\frac{1}{2}(\alpha+\beta)I_n+B$, we have
	$$
	B = \left[\begin{matrix}
	\,w I_{r}  & C\\
	D & -w I_{n-r}\,\end{matrix}\right]
	$$
and	we obtain  
    $$ H_{\theta}(B) \,=\, \left[\begin{matrix}
		   \,w_\theta I_r & N_\theta\\
			N_\theta^* & -w_\theta I_{n-r}\,\end{matrix}\right],	
			$$
	with $N_\theta=\frac{1}{2}\left(e^{-i\theta}C+e^{i\theta}\,D^*\right)$.
	By the singular value decomposition, there exist unitary matrices $U\in M_r$, $V\in M_{n-r}$, such that
	$$
	U^*\,N_\theta\, V\,=\, 	D_\theta\,\in\, M_{r,n-r}
	$$ 
	contains a diagonal principal submatrix with  the nonzero singular values    $s_1(\theta), \dots, s_p(\theta)$  of  $N_\theta$ in its main diagonal.
	Since the matrix $W=U\oplus V$ is unitary and  $W^* = U^*\oplus V^*$, we easily get
 $$
 		W^* H_{\theta}(B)\,W \, =\, \left[\begin{matrix}
		     w_\theta I_r & U^*N_\theta\,V\\
			\,V^*N_\theta^*\,U & -w_\theta I_{n-r}\,
			\end{matrix}\right]
		\,=  \,\left[\begin{matrix}\, w_\theta I_{r}& D_\theta\\
			D_\theta & - w_\theta I_{n-r}\,\\
		\end{matrix}\right].
 $$
 Thus, $W^* H_{\theta}(B)\,W$ is permutationaly similar to the direct sum of the $p$ blocks 
 $B_j \in M_2$, 	plus $r-p$ blocks $w_\theta$ and $n-r-p$ blocks $-w_\theta$ of  size $1$. 
	\smallskip
	
\noindent {\bf (iii)}  
 If $n<2r$ ($n>2r$), then it follows from (ii) that the eigenvalues of $H_\theta (B)$ are $w_\theta$  (resp.\! $-w_\theta$) and
\begin{equation}\label{vpH_B}
 \pm\sqrt{w_{\theta}^2+s_j^2(\theta)}, \qquad j=1, \ldots,n-r ,  
\end{equation}
(resp.\! $j=1, \ldots,r$). If $n=2r$, they are those in (\ref{vpH_B}). 
Having in mind (i) and  that 
\begin{equation}\label{HthetaAB}
	H_\theta(A)\,=\,\tfrac{1}{2}\Re(e^{-i\theta}(\alpha+\beta))I_n+H_\theta(B),
\end{equation}
the claimed eigenvalues of $H_\theta(A)$  are readily obtained.	
\end{pf} 

\medskip
\smallskip

Some considerations are in order. 
If   $\Re(e^{-i\theta}\alpha)$ or $\Re(e^{-i\theta}\beta)$ are eigenvalues of $H_\theta(A)$,
then they correspond to tangent lines of $C(A)$  passing through the point $\alpha$ or $\beta$, respectively, this meaning that these points 
belong to $C(A).$ Moreover, the remaining tangent lines form a family central symmetric relatively to $(\alpha+\beta)/2$ 
as implied by (\ref{lambdaj}).

The following technical lemma is used in the proof of the next theorem.

\medskip
\smallskip

\begin{lem}\label{discel}  Let $A\in M_n$ with two of the  eigenvalues of
	$H_\theta(A)$ of the form
	$$\Re (e^{-i\theta}c)\pm \tfrac{1}{2} \sqrt{a^2-(a^2-b^2)\sin^2(\theta- \phi)}, \qquad \theta\in[0,2\pi),$$
	for some $a,b,\phi \in \mathbb{R}$ and $c\in \mathbb{C}$. Then $C(A)$ contains an ellipse centered at $c$, with the major axis parallel to the vector $e^{i\phi}$ of length $|a|$ and the minor axis of length $|b|$. 
	\end{lem}

\medskip

If the non-diagonal blocks $C, D$ in \eqref{bloco} are such that $Z=DC$ is normal and commutes with $H=C^*C+DD^*$, then these matrices can be simultaneously diagonalized by the same unitary matrix, and we denote their eigenvalues by
	$z_j$ and $h_j$, $j=1, \dots, n-r$, labeled
	according to the order in which they appear in the respective diagonal matrices unitarily similar to $Z$ and $H$.
	
	 The next result is a refinement of \cite[Theorem 3.1]{GS2021} and is instrumental for our purposes. The proof takes the original ideas there combined with the introduced refinement in Lemma \ref{eigs}, ensuring the unitary reducibility  of $\Re(e^{-i\theta}A)$ to a direct sum of small sized blocks.

As usual, ${\rm Arg}\, z$  denotes the principal argument of the complex number $z$. 
	
	\medskip
	
\begin{thm}\label{principal}  Let $A\in M_n$ be of type \eqref{bloco}, $n\leq 2r$, such that 
 $Z=DC$ is normal and commutes with $H=C^*C+DD^*$, the respective eigenvalues  $z_j$ and $h_j$,  $j=1,\ldots,n-r$, being labeled as mentioned above.  Let
$\hat{\mathcal{E}}_j$ be the ellipse  with foci at 
$\frac{1}{2}(\alpha + \beta)\pm \frac{1}{2}\sqrt{\Delta_j}$,
for $\Delta_j=(\alpha-\beta)^2+4z_j$, 
with major and minor axes of length
	\begin{equation}\label{eixos}
	\sqrt{\tfrac{1}{2}\,|\alpha-\beta|^2+h_j+\tfrac{1}{2}\,|\Delta_j|}
	\qquad \hbox{and} 	\qquad 
	\sqrt{\tfrac{1}{2}\,|\alpha-\beta|^2+h_j-\tfrac{1}{2}\,|\Delta_j|}\,,
	\end{equation}
respectively. Then the following statements hold.
\begin{itemize}
\item[\bf (a)] The boundary generating curve of $W(A)$ is given by
$$C(A)\,=\,
	\left\{
	\begin{array}{ll} 
		\hat{\mathcal{E}}_1\cup\ldots\cup\hat{\mathcal{E}}_{n-r}\cup\{\alpha\}\,, &\mbox{if } \ n< 2r\\
		\hat{\mathcal{E}}_1\cup\ldots\cup\hat{\mathcal{E}}_{n/2}\,, & \mbox{if } \ n=2r\\
	\end{array}
	\right.$$
	and  $W(A)$ is the convex hull of $\hat{\mathcal{E}}_1,\dots,\hat{\mathcal{E}}_{n-r}$.
\medskip
 	
\item[\bf (b)] 
The higher rank numerical range of $A$ is given by
\begin{equation}\label{notsingleempty}
 \Lambda_k(A)\,=	\bigcap_{1 \leq j_1 \leq \cdots\leq j_{n-r-k+1} \leq n-r} \conv\, \{ \hat{\mathcal{E}}_{j_1}, \dots, \hat{\mathcal{E}}_{j_{n-r-k+1}} \}\,, \ \ \quad 1 \leq k\leq n-r,
\end{equation}
 and 
\begin{equation}\label{singleempty}
\hbox{either} \ \  \Lambda_k(A) =
	\left\{
	\begin{array}{cl} 
        \!\{\alpha\}\,,   &  n-r < k   \leq  r\\
        \emptyset\,,   &  r<k\leq n < 2r
	\end{array}
	\right.
\ \ \hbox{or} \ \quad \Lambda_k(A) =\,\emptyset\,, \ \ \frac{n}{2} <k\leq n=2r\,.
\end{equation}
\end{itemize}
\end{thm}

\begin{pf} {\rm (a)}  Since  $DC$ is normal and commutes with
	$C^*C+DD^*$, they can be diagonalized by the same unitary similarity, as already mentioned. 
	Under the theorem hypothesis,  if $n\leq 2r$, then the matrix
	$M_{C,D}(\theta)$
	is unitary diagonalizable for all values of $\theta$ by the same unitary similarity  of $DC$ and $C^*C+DD^*$, and so
	 its eigenvalues  are
\begin{equation}\label{eq1}
	\mu_j(\theta)\, 
	= \, h_j +2\,\Re(z_j)\cos(2\theta) +2\,\Im(z_j)\sin(2\theta), \qquad j=1, \dots, n-r.
\end{equation} 
	For $\theta \in [0, 2\pi)$, recalling   \eqref{w} and using 
$$
 2\,(\Re(w))^2=\Re(w^2)+|w|^2,\ \ \quad  2\,(\Im(w))^2=-\Re(w^2)+|w|^2,\ \ \quad 2\,\Re(w)\Im(w)=\Im(w^2),
$$
 by simple computations we get
\begin{eqnarray*}
	4\,w_\theta^2 +\mu_j(\theta) 
		&=& 4\,\big(\Re(w)\cos\theta+\Im(w)\sin\theta\big)^2 +\mu_j(\theta) \\
		&=& 4\,(\Re(w))^2\cos^2\theta+4\,(\Im(w))^2\sin^2\theta+4\,\Re(w)\Im(w)\sin(2\theta) 
		 +\mu_j(\theta) \\
	 & =&  2\,|w|^2+2\,\Re(w^2)\cos(2\theta)+2\,\Im(w^2)\sin(2\theta) +\mu_j(\theta)\\
	&=&  2\,|w|^2+h_j+2\,\Re(w^2+z_j)\cos(2\theta)+2\,\Im(w^2+z_j)\sin(2\theta) \\
	&=&  2\,|w|^2+h_j+2\,|w^2+z_j|\cos(2\theta-{\rm Arg}(w^2+z_j))\\
	&=&  2\,|w|^2+h_j+2\,|w^2 +z_j|-4\,|w^2+z_j|\sin^2(\theta-\phi_j)\,,
\end{eqnarray*}
with $2\phi_j={\rm Arg}(w^2+z_j)$.
Recalling  the expressions of the eigenvalues  of $H_\theta (A)$ in Lemma~2.1~(iii) for $\theta \in[0, 2\pi)$ and by Lemma~\ref{discel}, we easily conclude that,
for any $j\in \{1,\ldots,n-r\}$, to the pair of eigenvalues in \eqref{eigs}  corresponds  a component of $C(A)$, namely the  ellipse $\hat{\mathcal{E}}_j$ centered at $\frac{1}{2}(\alpha +\beta)$,  with major and minor axes of lengths
$$
 \sqrt{2\,|w|^2+h_j+2\,|w^2 +z_j|} \qquad \hbox{and} \qquad \sqrt{2\,|w|^2+h_j+2\,|w^2 +z_j|},
$$
respectively, and whose major axis is parallel to $e^{i \phi_j}$, this implying that $\hat{\mathcal{E}}_j$ has the asserted foci. 
If $n<2r$, then $\Re(e^{-i\theta}\alpha)$ is an eigenvalue of $H_\theta(A)$, and so the point $\alpha$ is also in $C(A)$.
Therefore, the statement  on $C(A)$ holds.
We easily see that $\alpha$ belongs to the elliptical discs bounded by $\hat{\mathcal{E}_j}$, $j=1, \dots, n-r$, and so, by Kippenhahn result, we have 
$$
  W(A)=\conv\, C(A)=\conv\,\{\hat{\mathcal{E}}_1, \dots, \hat{\mathcal{E}}_{n-r}\}.
$$ 
	
(b) If $k=1$, the result is clear. Otherwise, the characterization of the higher rank numerical range in \eqref{i} is equivalent to
$$
 \Lambda_k(A) \,= \bigcap_{\theta\in [0, \pi]} \left\{z\in W(A): \  \lambda_{n-k+1}{(\theta)} \leq \Re(e^{-i\theta} z )\leq \lambda_k{(\theta)}\right\}.
$$
By  Lemma~\ref{eigs}~(iii),
if $k \leq n-r$, then $\lambda_k(\theta)$ is the $k$-th largest number in the set 
$$	
 \left\{ \tfrac{1}{2}{\Re \big(e^{-i\theta}}(\alpha+\beta)\big) +\tfrac{1}{2}\sqrt{4\,w_\theta^2+\mu_j(\theta)}: \ j=1, \dots, n-r\right\}
$$
and $\lambda_{n-k+1}(\theta)$ is the $k$-th smallest number in the set
$$	
 \left\{ \tfrac{1}{2}{\Re \big(e^{-i\theta}}(\alpha+\beta)\big) -\tfrac{1}{2}\sqrt{4\,w_\theta^2+\mu_j(\theta)}: \ j=1, \dots, n-r\right\}.
$$
Moreover, if  $n-r< k\leq r$, then  $\lambda_k(\theta)=\lambda_{n-k+1}(\theta)=\Re(e^{-i\theta}\alpha)$, and if $r< k\leq n$, then
$\lambda_k(\theta)<\lambda_{n-k+1}(\theta)$.
Thus,  \eqref{singleempty} holds. From the proof in (a), we conclude that  
\begin{equation} \label{suplines}
 \Re(e^{-i\theta} z)=\lambda_k(\theta) \qquad  \hbox{and}\qquad  \Re(e^{-i\theta} z)= \lambda_{n-k+1}(\theta)
\end{equation} 
are the support lines of the elliptical component $\hat{\mathcal{E}}_k$ in $C(A)$ perpendicular to the direction $\theta$, when $ k \leq n-r$,
and  the smallest higher rank numerical range of $A$ with nonempty interior is
clearly 
$$
 \Lambda_{n-r}(A)\,=\,\bigcap_{j=1}^{n-r}\, {\rm conv}\, \hat{\mathcal{E}}_j.
$$  
If $1<k<n-r$, taking all the possible convex hulls 
$$ 
 \conv\, \{ \hat{\mathcal{E}}_{j_1}, \dots, \hat{\mathcal{E}}_{j_{n-r-k+1}}\}, \qquad 1 \leq j_1 \leq \cdots \leq j_{n-r-k+1} \leq n-r,
$$
and then intersecting them is equivalent to consider all the elements $z$ in $W(A)$ in the stripes defined by the lines (\ref{suplines}), for all $\theta \in [0, 2\pi)$. This means that (\ref{notsingleempty}) holds too. 
\end{pf}

\medskip

As previously noticed, the statements of Theorem \ref{principal},  and its consequences, 
may be formulated for $n>2r$, interchanging $C,D,$ as well as $\alpha,\beta$ and $r,n-r$.

 We remark that Theorem \ref{principal} holds, in particular, if the non-diagonal blocks of $A$ in \eqref{bloco}  are such that $CD$ and $DC$ are both normal matrices.

\medskip
  
\section{Matrices with elliptical higher rank numerical range}\label{Sec3}

\medskip

In  this section,  classes of block matrices of the form \eqref{bloco}, yielding   elliptical higher rank numerical ranges 
are presented.
 When the discs bounded by the ellipses in $C(A)$ of  Theorem \ref{principal} (a) form a nested chain, we get the following generic corollary.

  \medskip
  
\begin{corol}\label{ellip}
Let $A$ be under the hypothesis of  Theorem \ref{principal} and let $\mathcal{E}_{j}$ be the closed elliptical disc bounded by the  ellipse  $\hat{\mathcal{E}}_{j}$ there  described. If
$\mathcal{E}_{n-r} \subseteq \cdots \subseteq \mathcal{E}_{1},$ 
then 
  $$\Lambda_k(A)={\mathcal{E}}_k, \qquad  1 \leq k\leq n-r.$$ 
\end{corol}

\begin{pf} Let $1 \leq k\leq n-r$. Since by  hypothesis the elliptical discs are nested, we have
 $$\conv\, \{ \hat{\mathcal{E}}_{j_1}, \dots, \hat{\mathcal{E}}_{j_{n-r-k+1}} \}=\mathcal{E}_{j_1}, \qquad 1 \leq j_1 \leq \cdots\leq j_{n-r-k+1} \leq n-r.$$ 
From  Theorem \ref{principal} (b), it follows that the intersection in (\ref{notsingleempty}) reduces to
$$
 \Lambda_k(A)\,=\,\bigcap_{j_1=1}^k \mathcal{E}_{j_1} \, = \, \mathcal{E}_{k}.
$$
 \end{pf}

\begin{rmk}\label{Rmk} Recalling the proof of Theorem \ref{principal} (a) and the role of the eigenvalues of  $H_\theta(A)$ given in Lemma \ref{eigs}
(iii), 
it can be easily  seen that Corollary \ref{ellip} holds, 
whenever the eigenvalues of $M(\theta)$ in \eqref{eq1} satisfy $\mu_1(\theta)\geq \cdots\geq \mu_{n-r}(\theta)$  for  all $\theta \in [0, 2\pi)$. 
\end{rmk}

\medskip

If  $DC$ is a scalar multiple of the identity,  the following result holds.
	
\medskip

\begin{corol}\label{teor DC escalar}  Let $A$ be of type \eqref{bloco} with   $DC = z_1 I_{n-r}$
 and $h_1 \geq \dots \geq h_{n-r}$ be the eigenvalues, counting the multiplicities, of $C^*C+DD^*$. 
Let $\mathcal{E}^1_k$ be the  elliptical disc with foci at 
$
\frac{1}{2}(\alpha+ \beta)\pm \frac{1}{2}\sqrt{\Delta}, 
$
for $\Delta = (\alpha-\beta)^2 + 4z_1$,
 major and minor axes of lengths
$$
 \sqrt{\tfrac{1}{2}|\alpha-\beta|^2+h_k+\tfrac{1}{2}|\Delta|} \qquad \hbox{and}\qquad \sqrt{\tfrac{1}{2}|\alpha-\beta|^2+h_k-\tfrac{1}{2}|\Delta|}.
$$
If $n< 2r$, then 
	$$\Lambda_k(A)=
	\left\{
	\begin{array}{ccl} 
		\mathcal{E}^1_k\,, & \hbox{if}  & 1\leq k\leq n-r\\
        \{\alpha\}\,,  & \hbox{if}  &  n-r < k   \leq  r\\
        \emptyset\,, & \hbox{if}  &  r<k\leq n
	\end{array}
	\right..$$ 
If $n=2r$, then 
	$$\Lambda_k(A)=
	\left\{
	\begin{array}{cll} 
		\mathcal{E}^1_k\, , & \mbox{if} &   1\leq  k \leq n/2 \\
		\emptyset\, , & \mbox{if} & n/2 < k \leq n
	\end{array}
	\right..
	$$ 
\end{corol}

\begin{pf}  By hypothesis, $z_1$ is the unique eigenvalue of $DC$ and  Theorem \ref{principal} (a)  characterizes $C(A)$. All the elliptical components of $C(A)$, which are the boundaries of $\mathcal{E}^1_j$, $j=1, \dots, n-r$, have the same  foci. 
Since $h_1\geq \dots \geq h_{n-r}$, the chain of inclusions $$\{\alpha\} \subseteq \mathcal{E}^1_{n-r}\subseteq  \cdots\subseteq \mathcal{E}^1_1$$ 
is ensured by the axes length of the ellipses.
 By 	Corollary~\ref{ellip}, we get $\Lambda_k(A)=\mathcal{E}^1_k$, for $1\leq k\leq n-r$. By Theorem \ref{principal} (b),  the cases of $\Lambda_k(A)$ being a singleton or an emptyset are obtained.
\end{pf}

\smallskip

If $\Delta =0$ in Corollary \ref{teor DC escalar},  
then $\Lambda_k(A)$ is the circular disc centered at $\frac{1}{2}(\alpha+\beta)$ with radius of length
\[
 \sqrt{\tfrac{1}{8}|\alpha-\beta|^2+\tfrac{1}{4}h_k} \qquad \hbox{if} \ \ 1\leq k \leq n-r.    
 \] 

For an arrowhead matrix of the form
\[ 
\left[
\begin{array}{ccc|c}
\!\alpha &  &  & c_1 \\
    & \! \ddots &  & \vdots \\
     &    &   \! \alpha & c_{n-1} \\ \hline
d_1  &  \cdots & \! d_{n-1} & \beta
\end{array}
\right] \qquad \hbox{or}\qquad 
\left[
\begin{array}{c|ccc}
\!\alpha & c_1 & \dots & c_{n-1} \\ \hline
     d_1 & \beta &  &  \\
     \vdots  &  &  \ddots  &  \\ 
     d_{n-1} &  & & \beta
\end{array}
\right]
\]
with zeros at the omitted  entries, the next case is  immediate (cf.\! \cite{CH2012}).
\medskip

\begin{corol} Let $A$ be an arrowhead matrix of type \eqref{bloco}, $n\geq 3$,  with  
$C={\bf c}$ and $D={\bf d}^T$  (or $C={\bf c}^T$ and $D={\bf d}$) for ${\bf c}, {\bf d} \in \mathbb{C}^{n-1}$.
Then 
	$$\Lambda_k(A)=
	\left\{
	\begin{array}{ccl} 
		\mathcal{E}\,, & \hbox{if}  & k=1\\
        \{\alpha\}  \  (\hbox{or} \ \{\beta\}) \,,  & \hbox{if}  &  2 \leq k   \leq  n-1\\
        \emptyset\,, & \hbox{if}  &  k=n
	\end{array}
	\right.\!,$$ 
where ${\mathcal{E}}$ is the elliptical disc with foci at 
$
\frac{1}{2}(\alpha+ \beta)\pm \frac{1}{2}\sqrt{\Delta}, 
$
for $\Delta = (\alpha-\beta)^2 + 4\,{\bf c}^T{\bf d}$,
major and minor axis of length
$$
\sqrt{\tfrac{1}{2}|\alpha-\beta|^2+\|c\|^2+\|d\|^2+\tfrac{1}{2}|\Delta|}
 \qquad \hbox{and}\qquad
\sqrt{\tfrac{1}{2}|\alpha-\beta|^2+\|c\|^2+\|d\|^2-\tfrac{1}{2}|\Delta|}.
$$\end{corol}

\begin{pf}
This is an immediate consequence of Corollary \ref{teor DC escalar} with $z_1={\bf c}^T{\bf d}$ and  $r=n-1$  if $C={\bf c}, D={\bf d}^T$ 
 (or $r=1$ if $C={\bf c}^T, D={\bf d}$), which trivially applies as $DC={\bf c}^T{\bf d}$ 
and $C^*C+DD^*=\|{\bf c}\|^2+\|{\bf d}\|^2$ (resp.\! $CD={\bf c}^T{\bf d}$ and $CC^*+D^*D=\|{\bf c}\|^2+\|{\bf d}\|^2$).
\end{pf}

\medskip

Let $[\alpha, \beta]$ denote the line segment joining the numbers $\alpha, \beta$ if $\alpha\neq \beta$, reduced to the singleton $\{\alpha\}$ if $\alpha=\beta$.

\medskip

\begin{corol}\label{corol D=kC*}  Let $A$ be of type \eqref{bloco} with $D=\zeta \, C^*$ for some $\zeta \in\C$. 
Let $s_1\geq\cdots\geq s_p$ be the non-zero singular values of $C$, $p=\rank (C)$, $m=\min\{n-r, r\}$ and $m'\!=\max\{n-r,r\}$.
Let ${\mathcal{E}}_k$ be the elliptical disc  with foci at 
$\frac{1}{2}(\alpha + \beta)\pm \frac{1}{2}\sqrt{\Delta_k}$,
for $\Delta_k=(\alpha-\beta)^2+4\zeta s^2_k$, major and minor axes of length 
$$
	\sqrt{\tfrac{1}{2}\,|\alpha-\beta|^2+(1+|\zeta|^2)s^2_k+\tfrac{1}{2}\,|\Delta_k|} \ \ \quad \hbox{and} \ \quad \sqrt{\tfrac{1}{2}\,|\alpha-\beta|^2+(1+|\zeta|^2)s^2_k-\tfrac{1}{2}\,|\Delta_k|}.
$$
If $n<2r$ ($n>2r$), then 
$$\Lambda_k(A)=
	\left\{
	\begin{array}{ccl} 
		{\mathcal{E}}_k & \hbox{if}  & 1 \leq k\leq p\\
        \hbox{$[\alpha, \beta]$}  & \hbox{if}  &  p < k \leq  m\\
        \{\alpha\} \  (\hbox{or} \ \{\beta\})  & \hbox{if}  &  m <  k  \leq  m'\\
        \emptyset & \hbox{if}  &  m' <k\leq n
	\end{array}
	\right..$$ 
If $n=2r$, then 
$$\Lambda_k(A)=
	\left\{
	\begin{array}{ccl} 
		{\mathcal{E}}_k & \hbox{if}  & 1 \leq k\leq p\\
       \hbox{$[\alpha, \beta]$}  & \hbox{if}  &  p < k \leq  n/2\\
        \emptyset & \hbox{if}  &  n/2 <k\leq n
	\end{array}
	\right..$$ 
        \end{corol}
\begin{pf} 
Under the hypothesis, for $n\leq 2r$, the matrices $Z,H$  defined in Theorem \ref{principal}  are
$$ 
Z=\zeta\,C^*C \qquad \hbox{and} \qquad H=(1+|\zeta|^2)\,C^*C.
$$ 
and the matrix in (\ref{Mtheta}) becomes a scalar multiple of $C^*C$.
Then  $C(A)$ is characterized in Theorem \ref{principal} (a), being
$\hat{\mathcal{E}}_1, \dots, \hat{\mathcal{E}}_p$ the non-degenerate ellipses  there described, with
$$ 
z_j=\zeta\,s_j^2, \qquad  h_j=(1+|\zeta|^2)s_j^2,\qquad  j=1,\ldots,p,
$$ 
 the non-zero eigenvalues of $Z,H$, respectively, and 
$\hat{\mathcal{E}}_{p+1}=\cdots= \hat{\mathcal{E}}_{n-r}=\{\alpha, \beta\}$.

Moreover, the  non-zero eigenvalues of $M(\theta)$ are  of the form 
$$
	\mu_j(\theta) =  \left((1+|\zeta|^2+2\,\Re(\zeta)\cos(2\theta)+2\,\Im(\zeta)\sin(2\theta)\right)s_j^2 , \qquad j=1, \dots, p,
$$	
where $s_1 \geq \cdots \geq s_p$ are the non-zero singulares values of $C$.
Then
$$\mu_1 (\theta) \geq\cdots\geq \mu_{p}(\theta) > \mu_{p+1}(\theta)=\cdots =\mu_{n-r}(\theta)=0, \qquad \theta \in [0, 2\pi).$$
By Remark \ref{Rmk},  the non-degenerate elliptical discs ${\mathcal{E}}_j$ bounded by  $\hat{\mathcal{E}}_j$, $j=1, \dots, p$, are nested and contain the line segment joining $\alpha$ and $\beta$:
$$\mathcal{E}_{n-r} =\cdots =\mathcal{E}_{p+1}=\,[\alpha, \beta]\, \subset\, \mathcal{E}_{p}\, \subseteq\, \cdots \,\subseteq\, \mathcal{E}_{1}.$$ 
 By Corollary \ref{ellip}, we conclude that
  $
  \Lambda_k(A)={\mathcal{E}}_k$, if  $1 \leq k\leq p$, and   $\Lambda_k(A)=[\alpha, \beta]$, if $p < k\leq n-r$.
Otherwise, the result follows by Theorem \ref{principal} (b).

 If $n>2r$,  interchanging $\alpha, \beta$, as well as $n-r,r$ and $C,C^*$ in the above proof, the result is obtained, because the singular values of $C$ and $C^*$ coincide.
\end{pf}

\medskip

The higher rank numerical range of tridiagonal matrices with special structure was characterized in  \cite{AAS2018}, generalizing  known results in the literature for the classical numerical range.

Let $T( {\bf d},{\bf a},{\bf c})$ denote a tridiagonal matrix, with main diagonal ${\bf a}$,  the first  upper descending diagonal  ${\bf c}$
and the first lower descending diagonal ${\bf d}$. 
Let  $B({\bf \tilde c})=T({\bf \tilde c_e},  {\bf \tilde c_o},0)$ be  the bidiagonal matrix asssociated to the vector  ${\bf \tilde c}=(\tilde  c_1, \tilde  c_2, \tilde  c_3,\tilde  c_4, \dots)$, where 
$${\bf \tilde c_o}=(\tilde c_1, \tilde c_3, \dots)\qquad \hbox{and}\qquad{\bf \tilde  c_e}=(\tilde  c_2, \tilde  c_4, \dots).$$ 

Corollary \ref{trid} provides a direct proof of the  numerical range result in \cite[Theorem~3.3]{BS2004} and also the characterization of the higher rank numerical range.

\medskip

\begin{corol}\label{trid}  Let $T=T({\bf d},{\bf a},{\bf c})\in M_n$  with ${\bf a}=(\alpha, \beta, \alpha,\beta, \ldots)$,
 ${\bf c}=(c_1, \dots, c_{n-1})$, ${\bf d}=(d_1, \dots, d_{n-1})$,  such that 
$$
 c_j=\zeta\, \overline{d_j},  \quad j\in J_1\subseteq \{1, \dots, n-1\} \ \quad \hbox{and} \ \quad   d_j=\zeta\, \overline{c_j}, \quad j \in J_2=\{1, \dots, n-1\}\setminus J_1\,, 
$$
 for some $\zeta \in \mathbb{C}$. For the vector $\bf \tilde c$ with components
 $$\tilde c_j=d_j, \ \ j\in J_1 \qquad \hbox{and} \qquad \tilde c_j=c_j, \ \ j\in J_2\,,$$
let the elliptical discs
 ${\mathcal{E}}_1, \dots, {\mathcal{E}}_p$ be  described as in Corollary \ref{corol D=kC*}, with $s_1\geq\cdots\geq s_p$ the non-zero singular values of the bidiagonal matrix $B({\bf \tilde c})$ and $p$ its rank.
 If $n=2r+1,$ then 
$$
 \Lambda_k(T)=
	\left\{
	\begin{array}{ccl} 
		{\mathcal{E}}_k\,, & \hbox{if}  & 1 \leq k\leq p\\
        \hbox{$[\alpha, \beta]$}\,,  & \hbox{if}  &  p < k \leq  r\\
        \{\alpha\}\,, & \hbox{if}  &  r <  k  \leq  r+1\\
        \emptyset\,, & \hbox{if}  &  r+1<k\leq n
	\end{array}
	\right.\!.
	$$ 
If $n=2r$, then 
$$
 \Lambda_k(T)=
	\left\{
	\begin{array}{ccl} 
		{\mathcal{E}}_k \,, & \hbox{if}  & 1 \leq k\leq p\\
       \hbox{$[\alpha, \beta]$} \,,  & \hbox{if}  &  p < k \leq  n/2\\
        \emptyset\,, & \hbox{if}  &  n/2 <k\leq n
	\end{array}
	\right.\!.
$$
\end{corol}

\begin{pf} 
Without loss of generality, we may suppose that  $J_1$ is the subset of  odd  numbers in $\{1, \dots, n-1\}$,
because interchanging any pair of corresponding off-diagonal entries of the tridiagonal matrix 
$T$, resulting into $\tilde T$, the higher rank numerical range remains unchanged, observing that $H_\theta(T)$ and $H_\theta(\tilde{T})$ are both tridiagonal matrices with the same spectra, for every $\theta \in [0,2\pi)$. 

Then the result follows  easily by Corollary \ref{corol D=kC*}, because $T$ is permutationaly similar, via the permutation matrix whose columns are the vectors of the canonical basis of $\mathbb{R}^n$ reordered into the odd indexed $e_1, e_3, \dots $, followed by the even indexed ones $e_2, e_4, \dots$, to
\[
\left[\begin{matrix}\alpha\, I_r &C\\
        \,\zeta\, C^*&\beta\, I_{n-r}\,
        \end{matrix}\right]
\]
with $r={\lceil \frac{n}{2} \rceil}$ 
and $C=B({\bf \tilde c})$. 
\end{pf}

\smallskip

Considering in Corollary \ref{trid} that ${\bf c}={\bf 1}$ is the  vector with all the entries equal to $1$ and $\zeta=-1$, then 
the result in \cite[Theorem 2]{CH2001} for the numerical range of a continuant matrix with biperiodic main diagonal is extended to the rank-$k$ numerical range, since the singular values of the bidiagonal matrix $B({\bf 1})$ of order $\lceil\frac{n}{2}\rceil \times \lfloor\frac{n}{2}\rfloor$  are
(see e.g.\! \cite{BLS_Hyp}):
$$
  s_k = 2\cos\frac{k\pi}{n+1}, \qquad k=1,\ldots,\left\lfloor\tfrac{n}{2}\right\rfloor.
$$

A  {\it  Toeplitz matrix} is the one with constant  entries along each descending diagonal.
A {\it tridiagonal $2$-Toeplitz matrix}   is of the form $T({\bf d}, {\bf a}, {\bf c})$ with   biperiodic descending diagonals ${\bf a}$, ${\bf c}$, ${\bf d}$. 
 If the matrix $T$ in Corollary \ref{trid} is a tridiagonal $2$-Toeplitz matrix, such that either $c_1\overline{c_2}=\overline{d_1}d_2$ or $\overline{c_1}d_2=\overline{c_2}d_1$, then we find \cite[Theorem 6 and 11]{AAS2018}, where the singular values $s_k$ of the corresponding bidiagonal matrices $B({\bf \tilde c})$ are explicitly given by 
$$
  s^2_k = |c_1|^2+|d_2|^2 + 2|c_1 d_2|\,\cos\frac{2k\pi}{n+1}, \qquad k=1,\ldots,\left\lfloor\tfrac{n}{2}\right\rfloor.
$$
 
Let $A({\bf d}, {\bf a},{\bf c})$ denote the antitridiagonal matrix, with main  antidiagonal ${\bf a}$,  first lower and upper ascending  diagonais  ${\bf c}$ and ${\bf d}$, respectively, and zeros elsewhere.
Antitridiagonal matrices  $A({\bf d}, {\bf a},{\bf c})$ with at most two non-zero antidiagonals, under the conditions of \cite[Theorem 2 or 3]{BF2021}, have also elliptical higher rank numerical ranges, since they are permutationally similar to tridiagonal matrices with zero main diagonal, satisfying the hypothesis of the last corollary, as a consequence of \cite[Theorems 1 and 2]{BF2020}.

The ({\it backward})  {\it shift}   operator on $\mathbb{C}^n$ is represented by the $n$-square matrix with ones  on the  subdiagonal (resp.\! superdiagonal) and zeros elsewhere.
Henceforth, as an obvious consequence of Corollary \ref{trid}, we get the following  result obtained in \cite{Shift2012, SIAM}.
\medskip

\begin{corol}
The rank-$k$ numerical range of the $n$-dimensional shift operator, $n \geq 2$, is the  circular disc centered at the origin  with radius $cos \frac{k\pi}{n+1}$, if $1\leq k \leq \lfloor\frac{n+1}{2}\rfloor$,
and  the emptyset, otherwise. 
\end{corol}

\smallskip

Next, we consider the case when the spectrum of $M(\theta)$  is independent of $\theta$.

Denote by $E(\alpha, \beta; s)$ the elliptical disc with foci at $\alpha$ and $\beta$,  minor axis of length  $s$ if $\alpha \neq \beta$, reduced to a circular disc centered at $\alpha$ of radius $s/2$ if $\alpha=\beta$. 

\medskip

\begin{thm}\label{Tindeptheta}  Let $A\in M_n$ be the block matrix in {\rm (\ref{bloco})}, such that the spectrum of $M(\theta)$ is independent of $\theta$. Let $s_1\geq \cdots \geq s_p$ be the non-zero singular values of $C+D^*$, counting multiplicities.
If $n<2r$ ($n>2r$), then 
$$\Lambda_k(A)=
	\left\{
	\begin{array}{ccl} 
		E(\alpha, \beta; s_k)\,, & \hbox{if}  & 1 \leq k\leq p\\
        \hbox{$[\alpha, \beta]$}\,,  & \hbox{if}  &  p < k \leq  m\\
        \{\alpha\} \,  \hbox{(or} \, \{\beta\}\hbox{)}\,, & \hbox{if}  &  m <  k  \leq  m'\\
        \emptyset\,, & \hbox{if}  &  m'<k\leq n
	\end{array}
	\right.\!,$$ 
where $m=\min\{n-r, r\}$ and $m'\!=\max\{n-r,r\}$. If $n=2r$, then 
$$\Lambda_k(A)=
	\left\{
	\begin{array}{ccl} 
		{E} (\alpha, \beta; s_k)\,, & \hbox{if}  & 1 \leq k\leq p\\
       \hbox{$[\alpha, \beta]$} \,,  & \hbox{if}  &  p < k \leq  n/2\\
        \emptyset\,, & \hbox{if}  &  n/2 <k\leq n
	\end{array}
	\right.\!.$$ 
In particular, $W(A)\,=\,E(\alpha, \beta;\|C+D^*\|)$.
\end{thm}

\begin{pf}  By hypothesis, the spectrum of  $M(\theta)$ is independent of $\theta$, thus it 
coincides with the spectrum of $M(0)$.  
We observe that \[M_{C,D}(0)\,=\, C^*C+DD^*+2\Re(DC)=(C+D^*)^*(C+D^*).\] 
Then the non-zero eigenvalues of
$M(0)$ are $s^2_1 \geq\cdots\geq s^2_p$.  
For $\theta \in [0, 2\pi)$, we have 
\[ 
(\Re((\alpha-\beta)\,e^{-i\theta}))^2 =\, |\alpha-\beta|^2\cos^2\!\left(\theta-\tau\right),
\]
with $\tau ={\rm Arg}(\alpha - \beta)$. As in Lemma \ref{eigs} (iii), the $k$-th largest eigenvalue of $H_\theta(A)$ is
$$
	\lambda_k(\theta) =
	\left\{
	\begin{array}{ccl} 
			\frac{1}{2}{\Re \big(e^{-i\theta}}(\alpha+\beta)\big) +\frac{1}{2}\sqrt{
|\alpha-\beta|^2\cos^2\!\left(\theta-\tau\right)+s^2_j}\, , & \mbox{if} &   1 \leq k \leq p\\
			\max\{ \Re (e^{-i\theta}\alpha ), \Re (e^{-i\theta}\beta)\} \, , & \mbox{if} &   p < k \leq m \\
			\Re (e^{-i\theta}\alpha ) \ \ \big(\hbox{or} \ \ \Re (e^{-i\theta}\beta)\big)\, , & \mbox{if} &   m < k \leq m' \\		
			\min\{ \Re (e^{-i\theta}\alpha ), \Re (e^{-i\theta}\beta)\} \, , & \mbox{if} &   m' < k \leq n-p \\
			\frac{1}{2}{\Re \big(e^{-i\theta}}(\alpha+\beta)\big) -\frac{1}{2}\sqrt{|\alpha-\beta|^2\cos^2\!\left(\theta-\tau\right)+s^2_{n-k+1}}\, , & \mbox{if} & n-p < k \leq n
	\end{array}
	\right.\!,
$$
when $n\leq 2r$ (resp.\! $n>2r$). 
By Lemma \ref{discel},  $C(A)$ contains (at most) $p$ non-degenerate (distinct) elliptical components, all centered at $\frac{1}{2}(\alpha+ \beta)$, with  major axis parallel to the vector $e^{i\tau}$ of length 
 \[\sqrt{|\alpha-\beta|^2+s^2_j}\]  
and minor axis of length $s_j$,  $j=1, \dots, p$.
As all the ellipses have the same foci at $\alpha,\beta$ and $s_1\geq \cdots\geq s_p$, 
by the characterization of the rank-$k$ numerical range  in \eqref{i}, 
we get 
$$\Lambda_k(A)\,=\,{E}(\alpha, \beta; s_k), \qquad 1\leq k\leq p.$$
  If  $M(0)$ is not full rank and $p <k\leq m$, then  the $k$-th and $(n-k+1)$-th eigenvalues of $H_\theta(A)$ are the maximum and  minimum of $\{\Re (e^{-i\theta}\alpha ),\Re (e^{-i\theta}\beta)\}$, respectively, yielding an additional  degenerated component   in $C(A)$ and $\Lambda_{k}(A)=[\alpha,\beta]$. 
  If $n<2r$ ($n>2r$),  it is clear that $\Lambda_k(A)=\{\alpha\}$ (resp. $\Lambda_k(A)=\{\beta\}$),  whenever $m < k \leq m'$.
It is obvious that $\Lambda_k(A)=\emptyset$, whenever  $m'<k\leq n$.
\end{pf}

\smallskip

\begin{rmk} 
By the proof of Theorem \ref{Tindeptheta},  $C(A)$ contains the boundaries of ${E}(\alpha, \beta; s_j)$, $j=1, \dots,\min\{r, n-r\}$ ($s_j=0$,  $j>p$)  and  additionally $\{\alpha\}$ if $n<2r$ (or $\{\beta\}$ if $n>2r$). Since the foci of the elliptical components in $C(A)$ are eigenvalues of $A$, the spectrum of $A$ reduces to  $\{\alpha, \beta\}$.
\end{rmk}

\medskip

{\it Quadratic matrices} are those with minimal polynomial of degree two, which are unitarily similar to matrices of type \eqref{bloco}, where $C$ can be chosen to be positive semidefinite and  $D = O$.
 Their higher rank numerical ranges are known \cite{Li2011, TsoWu1999} and these results follow from 
 Theorem \ref{Tindeptheta}, Corollary \ref{teor DC escalar}  or Corollary \ref{corol D=kC*}, being  the numerical range
equal to  $E(\alpha, \beta;  \|C\|).$

Theorem \ref{Tindeptheta} includes matrices not covered by Theorem \ref{principal} (see e.g.\! \cite[Example 3]{GS2021}).

Let $\nu\subseteq \{1, \dots, r\}$ and $\overline{\nu}=\{1, \dots, r\}\setminus \nu$.
We have $DC=O$, when the $i$-th rows  of $C$ are zero, for $i \in \nu$ and  the $j$-th  columns  of $D$ are zero for $j\in \bar{\nu}$.  
Thus, the results on the block matrices in \cite[Theorems 5--7]{CH2012} follow as  corollaries of  Theorem~\ref{Tindeptheta} or Corollary~\ref{teor DC escalar}, and  their elliptical higher rank numerical ranges can also be explicitly obtained, as in the result below.

\medskip

\begin{corol} \label{HRquadr}
Let $A$ be of type \eqref{bloco}  with zeros on the $i$-th rows (columns)  of $C$, for $i \in \nu$, and  the $j$-th columns (resp.\! rows) of $D$, for $j\in \overline{\nu}$, and let $p=\rank (C+D^*)$.   Then $\Lambda_k(A)$ is characterized  
as in Theorem \ref{Tindeptheta}.
 In particular, $$\Lambda_k(A)=E\,(\alpha, \beta;\sqrt{h_k}), \qquad 1\leq k\leq p,$$ 
 with  $h_1\geq \cdots\geq h_p$   the non-zero eigenvalues of $C^*C+DD^*$ (resp.\! $CC^*+D^*D$).
\end{corol}

\begin{pf} The statements follow readily 
from Theorem~\ref{Tindeptheta}, since the hypotheses imply that
 $M_{C,D}(\theta)$ (resp.\! $M_{D,C}(\theta)$) is independent of $\theta$ and  the non-zero singular values of $C+D^*$ are equal to  $\sqrt{h_k}$, $ k=1,\dots,p$.
\end{pf}

\medskip

As the next example shows, block matrices of type \eqref{bloco} can have elliptical higher rank numerical range   
even if  the blocks $C,D$ are not under the conditions of Theorem \ref{principal} and  even if the hypothesis  of Theorem \ref{Tindeptheta} is not satisfied.

\medskip

\begin{Exa}
The non empty higher rank numerical ranges of the matrix
\begin{equation}\label{ex}
A=\left[
\begin{array}{cccccc}
 2 & 0 & 0 & 2+2 i & 1-i & 0 \\
 0 & 2 & 0 & -i & -1+i & 0 \\
 0 & 0 & 2 & 0 & 0 & 4 \\
 \frac{i}{4} & 0 & 0 & 3 & 0 & 0 \\
 \frac{i}{4} & \frac{3}{4}+\frac{i}{4} & 0 & 0 & 3 & 0 \\
 0 & 0 & \frac{1}{16} & 0 & 0 & 3 \\
\end{array}
\right]
\end{equation}
form a nested chain of elliptical discs. 
 In Figure \ref{f1}, the boundaries of $\Lambda_k(A)$ and their foci are represented in black for  $k=3$, blue for $k=2$ and red for $k=1$. 

\medskip

\begin{figure}[!h]\label{f1}
\centering
\includegraphics[scale=0.7]{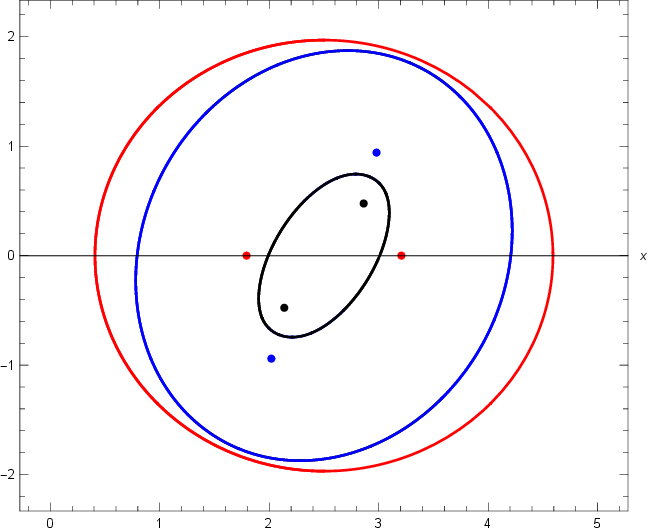}
\caption{Boundaries and foci of $\Lambda_3(A)\subset {\color{blue} \Lambda_2(A)} \subset {\color{red} \Lambda_1(A)}$ for $A$ in \eqref{ex}}
\end{figure}
\end{Exa}

Thus, the problem of finding necessary and sufficient conditions for occuring
elliptical shaped higher rank numerical ranges  for  block matrices of type \eqref{bloco}, that might include those of Theorem \ref{principal} and \ref{Tindeptheta}, can  still be raised.
Another question concerns the study of the higher rank numerical ranges  when the matrices \eqref{bloco}  are extended to bi-infinite  matrices  with two scalar infinite main diagonal blocks. 

\medskip
\medskip

\noindent {\bf Acknowledgments.} The authors thank the anonymous referee for carefully reading
the manuscript and helpful observations.

\medskip




\begin{thebibliography}{00}


  
\bibitem{AAS2018} M. Adam, A. Aretaki and I. M. Spitkovsky, Elliptical higher rank numerical range of some Toeplitz matrices, {\it Linear Algebra Appl.} {\bf 549} (2018), 256--275.
	
\bibitem{BdP1998} N. Bebiano and J. da Provid\^encia, Numerical ranges in physics, {\it Linear Multilin. Algebra} {\bf 43} (1998), 327--337. 
	
\bibitem{BPN2013} N. Bebiano, J. da Providência and A. Nata, The numerical range of banded biperiodic Toeplitz operators, 
 {\it J. Math. Anal. Appl.} {\bf 398} (2013), 189--198.

\bibitem{BLdP2004} N. Bebiano,  R. Lemos and J. da Providência, Numerical ranges of unbounded operators arising in quantum physics, {\it Linear Algebra Appl.} {\bf 381} (2004), 259--279.

\bibitem{BF2020} N. Bebiano and S. Furtado, 	Remarks on anti-tridiagonal matrices,
{\it Appl. Math. Comput.} {\bf 373} (2020), 125008.

\bibitem{BF2021} N. Bebiano and S. Furtado, A note on classes of structured matrices with elliptical type numerical range,
{\it Czechoslov. Math. J.} {\bf 71} (4) (2021),  1015--1023 
	
\bibitem{BPSK2021} N. Bebiano, J. da Providência, I. Spitkovsky and V. Kenya, Kippenhahn curves of some tridiagonal matrices, {\it Filomat} {\bf 35} (2021), 3047--3061.
	
\bibitem{BPS2022} N. Bebiano, J. da Providência and I. M. Spitkovsky, On Kippenhahn curves and higher-rank numerical ranges of some matrices,
{\it Linear Algebra Appl.} {\bf 629} (2021), 246--257.

\bibitem{BLS2023b} N. Bebiano, R. Lemos and G. Soares, On the numerical range of Kac-Sylvester matrices, {\it Electron. J. Linear Algebra} {\bf 39} (2023), 241--259.	
	
\bibitem{BLS_Hyp} N. Bebiano, R. Lemos and G. Soares, On the hyperbolicity of the Krein space numerical range, {\it Linear Multilin. Algebra} {\bf 72} (14) (2024), 2267--2287. 
		
\bibitem{BLS2024} N. Bebiano, R. Lemos and G. Soares, Algebraic curves associated with centrosymmetric matrices of orders up to $6$, {\it Adv. Oper. Theory} {\bf 9} (56) (2024).
	
\bibitem{BS2004} E. Brown and I. M. Spitkovsky,  On matrices with elliptical numerical ranges, {\it Linear Multilin. Algebra} {\bf 52} (2004), 177--193.
	
\bibitem{CKZ2006} M.-D. Choi, D. W. Kribs and K. \.{Z}yczkowski, Higher-rank numerical ranges and compression problems, {\it Linear Algebra Appl.} {\bf 418} (2006), 828--839.
	
\bibitem{CKZ2006A} M.-D. Choi, D. W. Kribs and K. \.{Z}yczkowski, Quantum error correcting codes from the compression formalism,
{\it Rep. Math. Phys.} {\bf 58} (1) (2006), 77--91. 
	 
\bibitem{CH2001}	 M.-T. Chien and J.-M. Huang, Numerical range of a continuant matrix,
{\it Applied Math. Lett.} {\bf 14} (2) (2001), 213--216,

\bibitem{CH2012} M.-T. Chien and K.-C. Hung, Elliptic numerical ranges of bordered matrices,  {\it Taiwanese J. Math.} {\bf 16} (3) (2012), 1007--1016.
	
\bibitem{CS2015} R. T. Chien and I. M. Spitkovsky, On the numerical ranges of some tridiagonal matrices,	{\it Linear Algebra Appl.} {\bf 470} (2015), 228--240. 
	
\bibitem{Shift2012} H. Gaaya, On the higher-rank numerical range of the shift, {\it J. Math. Sci. Adv. Appl.} {\bf 13} (2012), 1--19. 	
	
\bibitem{WG2013} H.~L.~Gau and P.~Y.~Wu, Higher-rank numerical ranges and Kippenhahn polynomials, {\it  Linear Algebra Appl.} {\bf 438} (2013), 3054--3061.
	
\bibitem{GS2021} T. Geryba and I.~M. Spitkovsky,  On the numerical range of some block matrices with scalar diagonal blocks, {\it Linear  Multilin. Algebra} {\bf  69} (2021), 772--785. 
	
\bibitem{GR} K. E. Gustafson and D. K. M. Rao, {\it  Numerical Range: The Field of Values of Linear Operators and Matrices}, Springer, 1997.
	
\bibitem{H1919} F. Hausdorff, Der Wertvorrat einer Bilinearform, {\it Math. Z.}  {\bf 3} (1919), 314--316.
		
	
\bibitem{JS2024} M. Jiang and I. M. Spitkovsky, Unified approach to reciprocal matrices with Kippenhahn curves containing elliptical components. {\it Linear Multilin. Algebra} {\bf 73} (7)  (2024), 1346--1368. 
	 
\bibitem{K1951} R. Kippenhahn, Über den Wertevorrat einer Matrix, {\it Math. Nachr.} {\bf 6}  (3-4) (1951),  193--228.

\bibitem{Li1996} C. K. Li,	A simple proof of the elliptical range theorem, {\it Proc. Amer. Math. Soc.} {\bf 124} (1996), 1985--1986.
	
\bibitem{Li2008} C.-K. Li and N.-S. Sze, Canonical forms, higher rank numerical ranges, totally isotropic subspaces, and matrix equations, {\it Proc. Amer. Math. Soc.}  {\bf 13} (2008), 3013--3023.
		
\bibitem{Li2009} C. K. Li,  Y.-T. Poon and N.-S. Sze, Condition for the higher rank numerical range to be non-empty. {\it Linear Multilin. Algebra}, {\bf 57} (4) (2009), 365--368.

\bibitem{Li2011} C. K. Li, Y.-T. Poon and N.-S. Sze, Elliptical range theorem for generalized numerical ranges of quadratic operators. {\it Rocky Mountain J. Math.}, {\bf 41} (3) (2011), 813--832.
	
\bibitem{JC2026} J. Ni\~{n}o-Cort\'{e}s and C. Vinzant, The convex algebraic geometry of higher-rank numerical ranges, {\it J. Symb. Comput.} {\bf 132} (2026), 102457.  
	
\bibitem{SIAM} E. Poon, I. M. Spitkovsky and H. J. Woerdeman, Factorization of singular matrix polynomials and matrices with circular higher rank numerical ranges, {\it SIAM J. Matrix Anal. Appl.}, {\bf 43} (3)
(2022), 1423--1439.

\bibitem{T1918}  O. Toeplitz,  Das algebraische Analogon zu einem Satze von Fej\'er, {\it Math. Z.} {\bf 2} (1918), 187--197. 		
		
\bibitem{TsoWu1999}  S--H. Tso and P. Wu, Matricial ranges of quadratic operators. {\it Rocky Mountain J. Math.}, {\bf 29} (3) (1999), 1139--1152.
	
\bibitem{W2008}  H. J. Woerdeman, The higher rank numerical range is convex. {\it Linear Multilin. Algebra}  {\bf 56} (1-2)  (2008), 65--67.  
	
\end{thebibliography}
\end{document}